\documentclass[10pt]{article}

\usepackage{myMacros}

\usepackage{amsmath, amssymb, amsthm, bm}
\usepackage[linesnumbered,ruled,vlined]{algorithm2e}
\usepackage{color}
\usepackage{graphicx}
\usepackage{subfig}

\begin{document}

\title{Surrogate Accelerated Bayesian Inversion for the Determination of the Thermal Diffusivity of a Material}

\author{
	James Rynn\footnote{Alan Turing Building, University of Manchester, Manchester, M13 9PL, UK and National Physical Laboratory, Teddington, TW11 0LW, UK (\texttt{james.rynn@manchester.ac.uk})},
	Simon Cotter\footnote{Alan Turing Building, University of Manchester, Manchester, M13 9PL, UK (\texttt{simon.cotter@manchester.ac.uk})},
	Catherine E Powell\footnote{Alan Turing Building, University of Manchester, Manchester, M13 9PL, UK (\texttt{catherine.powell@manchester.ac.uk})} 
	and Louise Wright\footnote{National Physical Laboratory, Teddington, TW11 0LW, UK (\texttt{louise.wright@npl.co.uk})}
}

\date{}

\maketitle

\begin{abstract}
Determination of the thermal properties of a material is an important task in many scientific and engineering applications.
How a material behaves when subjected to high or fluctuating temperatures can be critical to the safety and longevity of a system's essential components.
The laser flash experiment is a well-established technique for indirectly measuring the thermal diffusivity, and hence the thermal conductivity, of a material.
In previous works, optimization schemes have been used to find estimates of the thermal conductivity and other quantities of interest which best fit a given model to experimental data.
Adopting a Bayesian approach allows for prior beliefs about uncertain model inputs to be conditioned on experimental data to determine a posterior distribution,
but probing this distribution using sampling techniques such as Markov chain Monte Carlo methods can be incredibly computationally intensive.
This difficulty is especially true for forward models consisting of time-dependent partial differential equations.
We pose the problem of determining the thermal conductivity of a material via the laser flash experiment as a Bayesian inverse problem in which the laser intensity is also treated as uncertain.
We introduce a parametric surrogate model that takes the form of a stochastic Galerkin finite element approximation, also known as a generalized polynomial chaos expansion, and show how it can be used to sample efficiently from the approximate posterior distribution.
This approach gives access not only to the sought-after estimate of the thermal conductivity but also important information about its relationship to the laser intensity, and information for uncertainty quantification.
We also investigate the effects of the spatial profile of the laser on the estimated posterior distribution for the thermal conductivity. 
\end{abstract}

% Keywords:
%\noindent{\it Keywords\/}:
%Bayesian inverse problem,
%Markov chain Monte Carlo methods,
%measurement uncertainty,
%prior distributions,
%thermal diffusivity,
%stochastic Galerkin finite element methods,
%generalized polynomial chaos expansion

\section{Introduction}
\label{sec:intro}

Many measurements are made indirectly: quantities of interest (QoIs) are estimated by measuring one or more other quantities and calculating the required values from a model that links the QoIs to the measured quantities.
Examples include scatterometry, where surface profile is inferred from intensity measurements, determination of Young's modulus, where the measured quantities are force and displacement, and determination of thermal diffusivity from measurements of temperature and time. 

Thermal diffusivity is defined as 
\begin{equation}
	\alpha = \lambda / (\varrho \cp),
\label{eqn:a=l/pc}
\end{equation} 
where $\lambda$ is thermal conductivity, $\varrho$ is density, and $\cp$ is specific heat capacity.
Density and specific heat capacity can be measured independently, so measurement of thermal diffusivity is often used to evaluate thermal conductivity.
Thermal conductivity is a key property in understanding heat transport in solids.
Accurate characterisation of thermal conductivity supports design of thermal behaviour in products ranging from protective coatings for turbine blades in high temperature gas streams to  insulation for low-temperature carriers for transport of vaccines.
Reliable uncertainty evaluation of thermal conductivity enables designers to have confidence that their products will meet the required specification under all circumstances.

In some cases the models used in indirect measurements are simple and can be inverted to express QoIs explicitly in terms of the measured quantities,
but in many cases the models are not invertible (e.g. due to observational noise),
so the problem of determining the QoIs is more challenging and we have to solve an inverse problem.

Inverse problems present difficulties in applying the approach to uncertainty evaluation outlined in the GUM \cite{bipm2008evaluation} and its supplements.
The difficulties become more severe when the model used to link the measured quantities and the QoIs is computationally expensive, making sampling methods too time-consuming to use. One remedy is to replace the expensive model with a cheaper one that is sufficiently accurate for the task at hand, called a surrogate model \cite[Chapter 5]{rasmussen2015novel} (also known as an emulator or metamodel).
This paper reports on the application of an approach combining a parametric surrogate model with a Bayesian approach to uncertainty evaluation for the determination of thermal conductivity from the laser flash experiment.
This approach has made a computationally expensive problem tractable.

\begin{figure}
	\centering
	\includegraphics[width=0.725\textwidth]{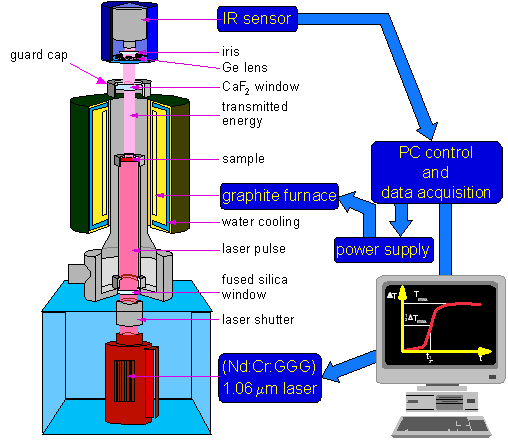}
	\caption{Diagram of the laser flash experiment set up.}
	\label{fig:experiment}
\end{figure}

In the rest of this section we describe the laser flash experiment (see Figure \ref{fig:experiment}),
present the mathematical model for the temperature of the material being tested and briefly discuss previous work by other authors.
We also outline our approach, linking it to existing literature and methods from other areas.
In Section \ref{sec:method}, we pose the problem of determining the thermal conductivity and the laser intensity as a Bayesian inverse problem.
We then describe how a Markov chain Monte Carlo (MCMC) algorithm that incorporates a surrogate forward model, here based on a stochastic Galerkin finite element method (SGFEM), can be used to efficiently probe the resulting approximate target distribution.
In Section \ref{sec:results}, we apply this methodology to a real life example where the sample material is a copper cylinder.
We discuss the construction of the SGFEM surrogate, timings, the estimated joint and marginal distributions for the thermal conductivity and laser intensity, and examine how these results can be interpreted with regards to the experimental data.
Finally, in Section \ref{sec:conclusions} we draw conclusions and briefly describe how this work may be extended to more complex problems.

We note that the laser flash experiment was also considered in a Bayesian setting by Allard et al in \cite{allard2015multi}.
Hence, since the set up of the physical experiment is similar, our discussion in Sections \ref{subsec:experiment} and \ref{subsec:mathmodel} below follows \cite[Section 1]{allard2015multi}.

\subsection{The laser flash experiment}
\label{subsec:experiment}

The thermal conductivity of a material is not a directly measurable quantity; its value is indirectly observed through the change in temperature of the material being characterised.
The laser flash experiment \cite{parker1961flash} is a commonly used method for making such observations.
In the experiment, the sample is placed in a furnace in a low-pressure inert gas atmosphere, supported by pins to minimise conductive heat losses, and heated to ambient temperature $\Ta$.
The lower face is then subjected to a short laser burst of duration $\tf$ seconds, heating the sample on that face.
The average temperature of the opposite face is measured by an infra-red sensor through a small window in the furnace.
The experimental set up is shown in Figure \ref{fig:experiment}.

\subsection{Mathematical model}
\label{subsec:mathmodel}

We assume that the sample is a cylinder $C$ of radius $R$ and vertical height $H$ and work in cylindrical coordinates $\rv=(r,\theta,z)$.
Furthermore, we assume the material is isotropic and homogeneous, with unknown scalar thermal conductivity $\lambda$ and known scalar density $\varrho$ and specific heat capacity $\cp$. 
Let $\zf$ denote the depth of penetration of the laser into the sample, which lasts for $\tf$ seconds.
Then the temperature of the material during the laser flash experiment is modelled by the solution $u$ to the transient heat equation
\begin{equation}
	\varrho \cp \partial_{t} u(\rv,t) = \nabla \cdot \left( \lambda \nabla u(\rv,t) \right) + Q(\rv,t), \qquad (\rv,t) \in C\times[0,T],
\label{eqn:pde}
\end{equation}
where the source term
\begin{equation}
	Q(\rv,t) := I \cdot \chi(r) \cdot 1_{\{[0,z_{f}]\times[0,\tf]\}}(z,t),
\label{eqn:q}
\end{equation}
represents the laser flash with unknown intensity $I$.
Here, the function $\chi\colon[0,R]\to\RR$ is either constant, to represent an ideal laser of uniform profile,
or of the form $\chi(r) = \exp(-r^{2}/2\rf^{2})$ to give a Gaussian profile in the radial direction as is the typical profile for lasers.

As stated in Section \ref{subsec:experiment}, the material is heated to the ambient temperature,
yielding the initial condition
\begin{equation}
	u(\rv,0) = \Ta, \qquad \rv \in C\cup\partial C.
\label{eqn:ic}
\end{equation}
We assume that the height $H$ is sufficiently small that heat losses across the curved surface are negligible.
Furthermore, we assume that heat is lost at the top and bottom faces at a rate proportional to the difference in temperature between the sample and the furnace.
This assumption is valid for convection and radiation at small temperature differences, as are generated experimentally.
Hence, we have the boundary conditions
\begin{align}
	\lambda \pdudn(\rv,t) & = 0,
	&& \rv \in \partial C_{R} := \left \{\rv \in \partial C \ \big| \ r=R\right\}, \ t \in [0,T], \label{eqn:bcR} \\
	\lambda \pdudn(\rv,t) & = \kappa \left( \Ta - u(\rv,t) \right),
	&& \rv \in \partial C_{H} := \left \{\rv \in \partial C \ \big| \ z\in\{0,H\}\right\}, \ t \in [0,T].
	\label{eqn:bcH}
\end{align}
We also impose the condition
\begin{equation}
	\lambda \pdudn(\rv,t) = 0,
	\quad
	\rv \in \partial C_{0} := \left \{\rv \in C \ \big| \ r=0\right\}, \ t \in [0,T]. \label{eqn:bc0}
\end{equation}
These conditions, together with the assumptions on the source term and material properties of the sample, ensure that the problem is axisymmetric.
That is, $u$ is symmetric about the $z$-axis and does not depend on $\theta$.

\subsection{Previous approaches}
\label{subsec:previousapproaches}

Many previous approaches to determining the unknowns $\lambda$ and $I$ have used optimization (e.g., see \cite{whittle1971optimization,wright2011parameter}). 
In such works, values which are by some measure optimal are found by matching model output to experimental data.
There are a number of drawbacks to such approaches.
Firstly, the optimal values that are returned do not have associated uncertainties.
The standard GUM \cite{bipm2008evaluation} approach to calculating the uncertainty is unsuitable because the problem is an inverse one, and the problem is too computationally expensive for standard Monte Carlo simulation.
Secondly, well-posedness of the problem in terms of uniqueness of the solution can be an issue, particularly in cases where there are many unknowns and there is little data to infer from.
A Bayesian formulation of the problem \cite{dashti2016bayesian,stuart2010inverse} allows both of these issues to be resolved.

Some previous work \cite{allard2015multi} has been done to treat the laser flash experiment in a Bayesian manner.
However, the posterior distribution for $\lambda$ and $I$ that is provided (up to a constant of proportionality) by the Bayesian approach must be evaluated for one to compute quantities of interest such as probabilities, expectations and variances.
Here, each evaluation of the posterior requires an evaluation of the solution $u$ of \eqref{eqn:pde}--\eqref{eqn:bc0}.
Since an exact solution is unavailable, $u$ must be approximated using a finite element method combined with a time-stepping scheme, a finite difference scheme or similar approach.
We may exploit axisymmetry to reduce the spatial domain to two dimensions.
However, depending on the choice of spatio-temporal discretization, and its level of fidelity, a single solve may still take tens of seconds, or minutes to perform on a standard desktop computer.
If computational resources are limited, then it is infeasible to perform a large number of forward solves.

Markov chain Monte Carlo (MCMC) methods \cite{brooks2011handbook} are popular algorithms for sampling from a distribution that is known only up to a constant of proportionality.
However, it is well documented that they suffer from a slow rate of convergence \cite{caflisch1998monte}, with the sampling error behaving as $\calO(M^{-1/2})$, where $M$ is the number of samples used.
The large number of samples required to estimate quantities of interest accurately, coupled with the expense of each posterior evaluation means that MCMC methods are often infeasible.
This is especially true when the forward problem consists of a time-dependent partial differential equation (PDE).

We will use a pre-computed parametric surrogate for the forward solution $u$ that can be cheaply evaluated for any choice of $\lambda$ and $I$.
This takes the form of a polynomial expansion and allows for fast sampling from the resulting approximate posterior distribution within an MCMC routine.
As well as reducing the time required to solve the inverse problem by several orders of magnitude, we find that with the proper choice of polynomial degree, there is no significant loss of accuracy compared to the standard approach using a fixed spatio-temporal discretization.

The idea of combining surrogate models with MCMC methods to accelerate the solution of Bayesian inverse problems is becoming popular in groundwater flow modelling \cite{domesova2017bayesian},  electrical impedance tomography \cite{hakula2014reconstruction,hyvonen2015stochastic,leinonen2014application} and other applications \cite{cui2015data,marzouk2009dimensionality,marzouk2007stochastic}.

\section{Method}
\label{sec:method}

In this section we discuss a parametric form of the forward problem \eqref{eqn:pde}--\eqref{eqn:bc0} and describe how its solution may be approximated using a stochastic Galerkin finite element method (SGFEM).
We also formally introduce the Bayesian inverse problem of determining $\lambda$ and $I$ from experimental data.
Finally, we explain how to combine the SGFEM surrogate with an MCMC method to approximate posterior distributions.

\subsection{Parametric forward problem and SGFEM}
\label{subsec:forwardproblem}

We consider the laser intensity $I$ as well as the thermal conductivity $\lambda$ to be unknown.
This is because the exact strength of the laser is not known and the effects of the laser depend on the heat absorption characteristics of the sample.
Note that we could treat other model inputs (such as density $\varrho$ or ambient temperature $\Ta$) as unknown and propagate all measurement uncertainties through the model but choose not to here for clarity.

To construct the surrogate, we start by modelling $\lambda$ and $I$ as real-valued random variables. Specifically, we assume that $\lambda$ and $I$ may be expressed in terms of independent uniform random variables with mean zero and unit variance.
That is,
\begin{equation}
	\lambda(\xi_{1}) = \mu_{\lambda} + \nu_{\lambda} \xi_{1}, \qquad I(\xi_{2}) = \mu_{I} + \nu_{I} \xi_{2},
\label{eqn:ylambdaI}
\end{equation}
for some $\mu_{\lambda}, \mu_{I}, \nu_{\lambda}, \nu_{I}\in\RR^{+}$ (to be chosen) with 
\begin{equation}
\xi_{1},\xi_{2} \sim \calU(-\sqrt{3},\sqrt{3}).
\end{equation}
For heterogeneous or layered materials, $\lambda$ may be more accurately represented as a spatial random field \cite{adler1981geometry}.
If the material consists of layers, each with a distinct unknown value of $\lambda$, then a random field could be constructed using a linear combination of spatial indicator functions with uncertain coefficients (as used for example in \cite{hyvonen2015stochastic}).
Alternatively, if the mean and covariance function of the random field are known, a truncated Karhunen-Lo{\`e}ve expansion \cite{loeve1978probability} could be used.
Although we do not consider $\lambda$ and $I$ to be spatially varying here,
the approach outlined below can be used whenever $\lambda$ and $I$ are linear functions of a finite set of independent random variables with appropriate distributions.

To set up the parametric forward problem, we make a change of variable.
Let $y_{1}$ and $y_{2}$ be the images of $\xi_{1}$ and $\xi_{2}$, respectively.
Then, $\yv:=(y_{1}, y_{2}) \in\Gamma$ where $\Gamma:=[-\sqrt{3},\sqrt{3}]^{2}$ is the parameter domain. Assuming that $\lambda$, $I$ are written as in \eqref{eqn:ylambdaI},
we may rewrite \eqref{eqn:pde}--\eqref{eqn:bc0} in the equivalent parametric form
\begin{align}
\begin{split}
\varrho \cp \partial_{t} u(\rv,t,\yv) & = \nabla \cdot \left( \lambda(y_{1}) \nabla u(\rv,t,\yv) \right) + Q(\rv,t,y_{2}), \\
& \hspace{3.8cm} (\rv,t,\yv) \in C \times [0,T] \times\Gamma,
\label{eqn:ypde}
\end{split}\\
u(\rv,0,\yv) & = \Ta,
\hspace{2.8cm} (\rv,\yv) \in \overbar{C}\times\Gamma,
\label{eqn:yic} \\
\lambda(y_{1}) \pdudn(\rv,t,\yv) & = \kappa \left( \Ta - u(\rv,t,\yv) \right),
\hspace{0.5cm} (\rv,t,\yv) \in \partial C_{H} \times [0,T] \times \Gamma,
\label{eqn:ybcH} \\
\lambda(y_{1}) \pdudn(\rv,t,\yv) & = 0,
\hspace{3cm} (\rv,t,\yv) \in \left( \partial C_{0} \cup \partial C_{R} \right) \times [0,T] \times \Gamma,
\label{eqn:ybc0}
\end{align}
where
\begin{equation}
	Q(\rv,t,y_{2}) := I(y_{2}) \cdot \chi(r) \cdot 1_{\{[0,z_{f}]\times[0,\tf]\}}(z,t).
\label{eqn:yq}
\end{equation}

We have chosen uniform random variables to ensure that weak formulations of the parametric model are well-posed.
We must also choose the values of $\mu_{\lambda}$ and $\nu_{\lambda}$ in \eqref{eqn:ylambdaI}  such that realisations of $\lambda$ are positive.
Note also that it is physically meaningless to model either $\lambda$ or $I$ as a random variable with non-zero probability of taking a negative value (e.g. Gaussian).

The solution $u$ to \eqref{eqn:ypde}--\eqref{eqn:yq} is a function of the parameters $y_{1}$ and $y_{2}$, as well as $\rv$ and $t$.
We will approximate it using a spectral stochastic finite element method \cite{babuska2004galerkin,ghanem2003stochastic} based on Galerkin approximation.
Specifically, we combine finite element approximation \cite{braess2007finite} on the physical domain with global polynomial approximation on the parameter domain.
For parabolic PDEs, the suitability of this approach has been analyzed in \cite{nobile2009analysis}.

If we exploit axisymmetry, then the SGFEM approximation takes the form 
\begin{equation}
	\uhk(\rv,t,\yv) := \sum\limits_{i=1}^{\nh} \sum\limits_{j=1}^{\nk} u_{ij}(t) \phi_{i}(\rv)\Psi_{j}(\yv),
\label{eqn:uhk}
\end{equation}
where $\phi_{i}$, $i=1,2,\dots,\nh$, are standard piecewise linear finite element functions associated with a mesh of triangular elements of width $h$ on the spatial domain
\begin{equation}
D:=\{\rv = (r,z) \ | \ r\in(0,R), \ z \in (0,H)\}.
\end{equation}
The functions $\Psi_{j}$, $j=1,2,\dots,\nk$, are chosen to form a basis for the set of polynomials in $y_{1}$ and $y_{2}$  on $\Gamma$ of total degree less than or equal to $k$.
Hence,
 \begin{equation}\label{nk_def}
 \nk = (k+2)(k+1)/2.
 \end{equation}
In particular, we use Legendre basis polynomials which are orthonormal with respect to the inner product
\begin{equation}
\langle \Psi_{j}, \Psi_{s} \rangle_{\rho} := \int_{\Gamma} \Psi_{j}(\yv) \Psi_{s}(\yv) \rho(\yv) \dyv,
\end{equation}
where  $\rho(\yv)=(2\sqrt{3})^{-2}$ is the joint probability density function of $\xi_{1}$ and $\xi_{2}$. We may then also write (\ref{eqn:uhk}) as 
\begin{equation}
	\uhk(\rv,t,\yv) := \sum\limits_{j=1}^{\nk} u_{j}(\rv,t) \Psi_{j}(\yv),
\label{eqn:uhk-B}
 \end{equation}
which is known as a polynomial chaos expansion \cite{xiu2002modeling}.

To find the functions $u_{ij}(t)$ in (\ref{eqn:uhk}) we solve the associated Galerkin equations \cite{benner2015low} which leads to a coupled system of ODEs.
For this, we use an implicit Euler scheme with uniform step size $\tau:=T/n_{t}$.
This yields a sequence of discrete approximations
\begin{equation}
	\uhkt^{n}(\rv,\yv) := \sum\limits_{i=1}^{\nh} \sum\limits_{j=1}^{\nk} u_{ij}^{n} \phi_{i}(\rv)\Psi_{j}(\yv),
\label{eqn:uhkt}
\end{equation}
at the time steps $\tau_{n}=n\tau,$ $n = 0,1, \dots,\nt$, and requires the solution of $\nt$ linear systems of dimension $n_{h}n_{k}$.
We will use the notation $\uhkt$ to denote a continuous function of time which interpolates $\uhkt^{n}$ at each $\tau_{n}$.

Once the coefficients $u_{ij}^{n}$ in \eqref{eqn:uhkt} have been computed, $\uhkt^{n}$ may be evaluated at any $\yv \in \Gamma$ of interest by evaluating the polynomials $\Psi_{j}$.

\subsection{Bayesian inverse problem formulation}
\label{subsec:bip}

Suppose we wish to recover the unknown thermal conductivity $\lambda$ given the (thermogram) results from a laser flash experiment with uncertain laser intensity.
We also need to recover the laser intensity $I$ since it directly affects the observed thermogram.
Specifically, we study a cylindrical sample of pure copper.
Let $t_{1},t_{2},\dots,t_{\nd}$ denote the \textit{measurement times} and $D_{L}$ denote a disc of radius $L$ centrally positioned on the top face of the sample over which an average temperature is measured at each $t_{n}$.
Our data $\dv\in\RR^{\nd}$ is the vector of observed averaged temperature values.
The data used here was measured at NPL in 2004 on a Netzsch LFA427 using a sample of copper at ambient temperature $\Ta = 385$\,K.

The Bayesian approach to this inverse problem is to condition our prior knowledge regarding $\lambda$ and $I$ on the measurements $\dv$.
We may obtain this knowledge from historical information, expert elicitation or even previous Bayesian analysis using another data set.
We assume that our measurements are subject to independent, identically distributed (iid), mean zero Gaussian noise, the variance $\sigma^{2}$ of which is also unknown.
In this section, we choose a prior distribution, formulate the Bayesian inverse problem and construct the posterior distribution for the unknowns.

\subsubsection{Prior distribution.}
\label{subsubsec:prior}

\begin{figure}
	\centering 
	\subfloat[]{\includegraphics[width=\textwidth]{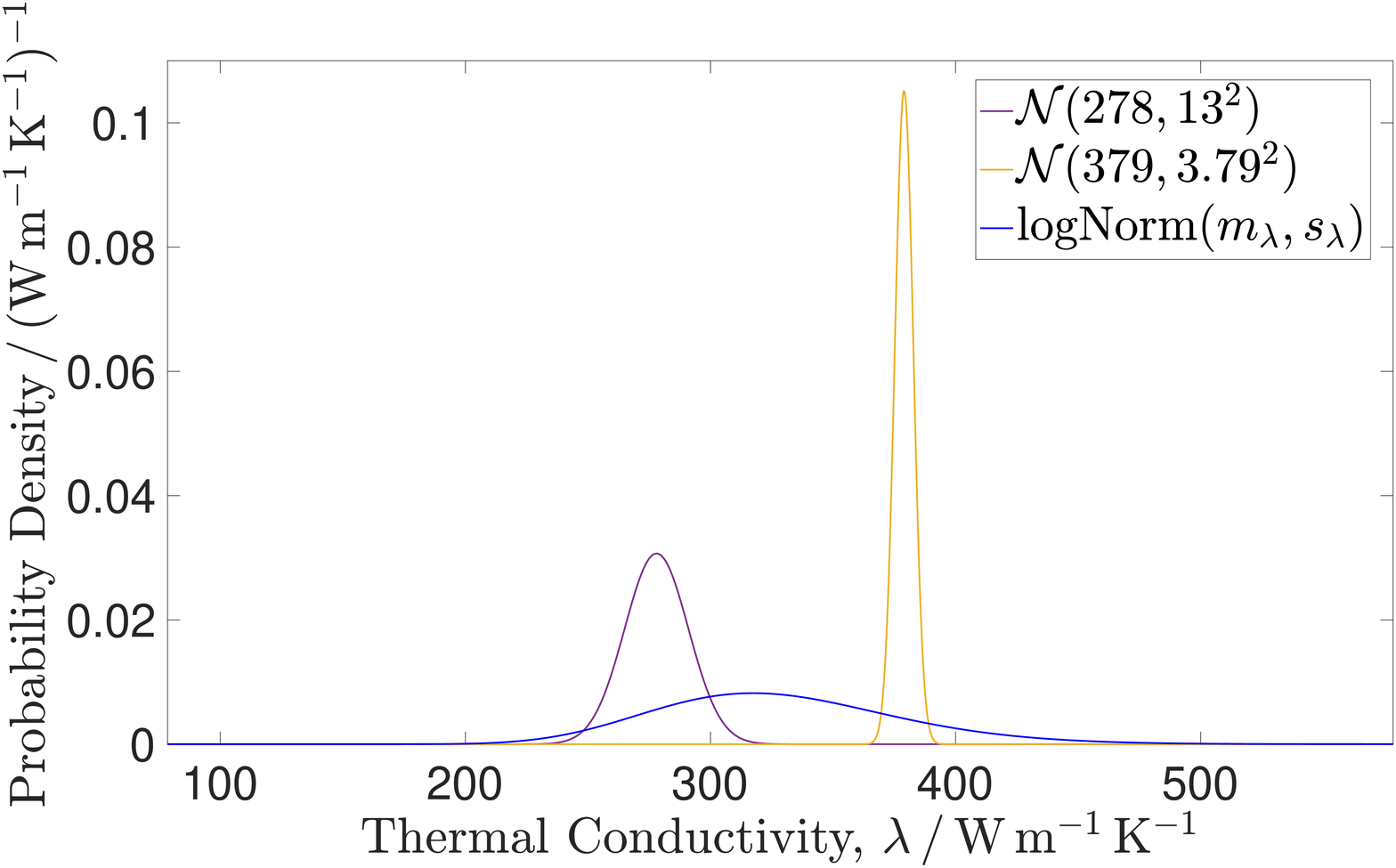}\label{fig:figure2a}} \\
	\subfloat[]{\includegraphics[width=\textwidth]{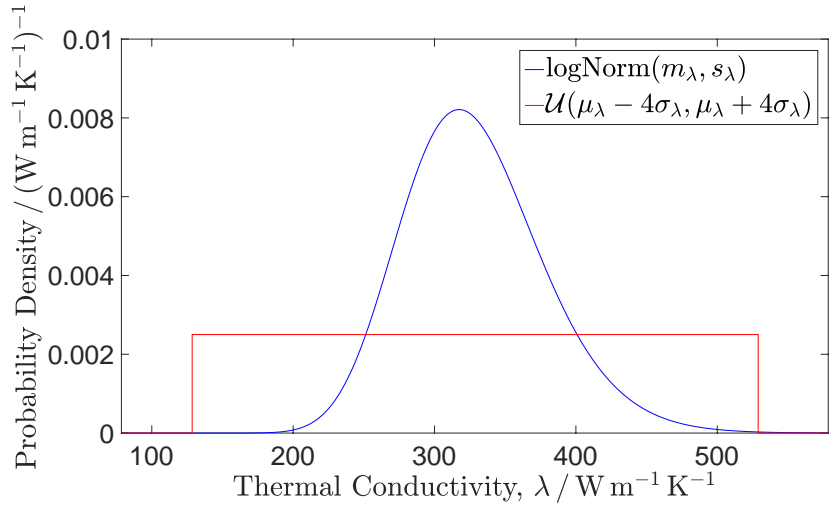}\label{fig:figure2b}}
	\caption{
		Distributions on the thermal conductivity $\lambda$:
		(a) two previously obtained Gaussian distributions and logNormal prior distribution.
		(b) logNormal prior distribution, and uniform distribution used to construct surrogate.}
\end{figure}

First, we assign a prior distribution to $\lambda$, which describes our belief about its value before incorporating the data $\dv$.
Two previous studies of samples involving copper deduced the distributions $\calN(379, 3.79^{2})$ and $\calN(278,13^{2})$.
In the first case, the value of 379 W\,m$^{-1}$\,K$^{-1}$ was obtained by directly fitting a model to data measured on a uniform sample of copper, and the estimated standard deviation of 3.79 W\,m$^{-1}$\,K$^{-1}$ was based on the typical uncertainty for the measurement.
In the second case, the value of 278 W\,m$^{-1}$\,K$^{-1}$ was obtained from analysis of a layered sample that was half braze and half copper, and a Monte Carlo analysis was used to propagate several other measurement uncertainties through a model to give the standard deviation of 13 W\,m$^{-1}$\,K$^{-1}$.
The second estimate is considered to be less reliable because it is affected by assumptions made about the properties of the braze.

As high probability confidence intervals from the two stated Gaussian distributions do not overlap (see Figure \ref{fig:figure2a}),
we have conflicting beliefs about the value of the thermal conductivity.
This leads us to be less certain about the value of the thermal conductivity than either of the previous studies suggests in isolation and makes construction of a prior both challenging and subjective.
To combine the two results into a prior distribution, in a way which reflects this uncertainty, we could consider a Gaussian distribution with (averaged) mean $\mu_{\lambda}:=(379+278)/2 = 328.5$ W\,m$^{-1}$\,K$^{-1}$ and large standard deviation $\sigma_{\lambda}:=50$ W\,m$^{-1}$\,K$^{-1}$.
However, to ensure positivity, we will use
\begin{equation}
	\lambda \sim \mbox{logNormal}(m_{\lambda}, s_{\lambda}),
\label{eqn:priorlambda}
\end{equation}
as our prior and perform inference on
\begin{equation}
	\theta_{1} := \ln(\lambda) \sim \calN(m_{\lambda},s_{\lambda}^{2}),
\label{eqn:theta1}
\end{equation}
where $m_{\lambda}, s_{\lambda}$ are chosen so that $\EE[\lambda]=\mu_{\lambda}$ and $\mbox{sd}[\lambda]=\sigma_{\lambda}$ under the prior.
That is, our prior density on $\theta_{1}$ is
\begin{equation}
	\pi_{0}^{\lambda}(\theta_{1}) = \frac{1}{\sqrt{2\pi s_{\lambda}^{2}}}\exp\left(-\frac{1}{2s_{\lambda}}(\theta_{1} - m_{\lambda})^{2}\right).
\label{eqn:prior_lam}
\end{equation}
The effects of the choice of prior are discussed later in this section.

As we have no information about the value of the laser intensity $I$ other than that it is positive, we choose the improper prior distribution ${\cal U}(0,\infty)$ with ``density'' $\pi_{0}^{I}(I) \propto 1_{[0,\infty)}(I)$.
For ease of sampling, we work with $\theta_{2} := \ln(I)$ and our prior ``density'' on $\theta_{2}$ is therefore given by
\begin{equation}
	\pi_{0}^{I}(\theta_{2}) \propto \exp(\theta_{2}).
\label{eqn:prior_I}
\end{equation}
This is an example of an uninformative prior.
Another common choice is the Jeffreys' prior \cite{jeffreys1998theory}, which is based on the Fischer information matrix \cite{allard2015multi,lehmann2006theory}.

We choose an inverse gamma prior distribution for $\sigma^{2}$ to ensure positivity and so that we can exploit it's conjugacy with the Gaussian distribution.
That is, our prior density on $\sigma^{2}$ is given by
\begin{equation}
	\pi_{0}^{\sigma}(\sigma^{2}) = \frac{\beta_{\sigma}^{\alpha_{\sigma}}}{\Gamma(\alpha_{\sigma})} (\sigma^{2})^{-\alpha_{\sigma}-1} \exp\left( -\frac{\beta_{\sigma}}{\sigma^{2}} \right),
\label{eqn:prior_sig2}
\end{equation}
where, here, $\Gamma(\cdot)$ is the gamma function.
The hyper parameters $\alpha_{\sigma} = 3$ and $\beta_{\sigma} = 0.0079$ are chosen so that all reasonable values of $\sigma^2$ have relatively high probability with respect to the prior distribution.

Assuming (a priori) independence of $\lambda$, $I$ and $\sigma^{2}$ (and therefore of $\theta_{1}$, $\theta_{2}$ and $\sigma^{2}$), letting $\thetav:=(\theta_{1},\theta_{2})$
we define our prior density on $(\thetav$, $\sigma^{2})$ to be the product of the three densities \eqref{eqn:prior_lam}, \eqref{eqn:prior_I} and \eqref{eqn:prior_sig2}.
That is,
\begin{equation}
	\pi_{0}(\thetav,\sigma^{2}) = \pi_{0}^{\lambda}(\theta_{1}) \pi_{0}^{I}(\theta_{2}) \pi_{0}^{\sigma}(\sigma^{2}).
\label{eqn:priortheta}
\end{equation}
We note that choosing a prior may be viewed as a form of regularization \cite{kaipio2006statistical}, allowing us to assign greater weight to values we deem more likely and less weight to those we believe less likely. Numerical investigations showed that for our problem the data $\dv$ is sufficiently informative about the unknowns that the choice of prior is not important.
That is, changing the prior distribution (including the  hyper parameters) had little effect on the estimated posterior. These studies are not presented in this paper.
However, experiments involving different choices of prior for similar inverse problems are given in \cite{allard2015multi,heidenreich2015statistical}.

\subsubsection{Data, likelihood, posterior and target densities.}
\label{subsubsec:datalikelihoodposterior}

We assume that our thermogram data $\dv$ are indirect observations of $\thetav$, subject to additive noise.
That is,
\begin{equation}
	\dv = \calvG(\thetav) + \etav, \quad \etav \sim \calN(\mathbf{0},\sigma^{2}I_{\nd}),
\label{eqn:data}
\end{equation}
where $\calvG\colon\RR^{2}\to\RR^{\nd}$ is the \textit{observation operator}, mapping values of $\thetav$ into data and $\etav$ is the noise vector.
A diagonal covariance matrix $\Sigma:=\sigma^{2}I_{\nd}$ is chosen as we assume individual measurements are subject to iid $\calN(0,\sigma^{2})$ noise.
This yields a likelihood function $L\colon\RR^{\nd}\times\RR^{2}\times\RR^{+}\to\RR$ (describing the probability of observing $\dv$ for a specific choice of $\thetav$ and $\sigma^{2}$) given by
\begin{equation}
	L(\dv\mid\thetav,\sigma^{2}) \propto (\sigma^2)^{-\nd/2} \exp(-\Phi(\thetav,\sigma^{2};\dv)),
\label{eqn:likelihood}
\end{equation}
where
\begin{equation}
	\Phi(\thetav,\sigma^{2};\dv) := \frac{1}{2}|\dv-\calvG(\thetav)|_{\Sigma}^{2} = \frac{1}{2\sigma^{2}}\|\dv-\calvG(\thetav)\|_{2}^{2},
\label{eqn:loglikelihood}
\end{equation}
is the negative log-likelihood or so-called \textit{potential}.
The observation operator $\calvG$ is given explicitly by
\begin{equation}
	\calvG(\thetav) = \left( \bar{u}(t_{1};\thetav), \bar{u}(t_{2};\thetav), \dots, \bar{u}(t_{\nd};\thetav)  \right)^{\top},
\label{eqn:G}
\end{equation}
where
\begin{equation}
	\bar{u}(t;\thetav):=\frac{1}{|D_{L}|}\int_{D_{L}}u(\rv,t;\thetav)\mbox{d}S_{z}(\rv)
\label{eqn:ubar}
\end{equation}
is the average temperature over the disc $D_{L}$ on the top face of the sample at time $t$ and $u(\rv,t;\thetav)$ denotes the solution to the deterministic forward problem (\ref{eqn:pde}--\ref{eqn:bc0}) for a \textit{fixed} choice of $\thetav$ where $(\lambda,I) = \exp(\thetav):=(\exp(\theta_{1}), \exp(\theta_{2}))$.

Note that \eqref{eqn:likelihood} is simply the probability density function of a multivariate Gaussian distribution.
Now, from Bayes' Theorem \cite{lee2012bayesian}, the posterior density $\pid\colon\RR^{2}\times\RR^{+}\times\RR^{\nd}\to\RR^{+}$ (describing the probability of obtaining values for $(\thetav, \sigma^{2})$ given the data $\dv$) satisfies
\begin{equation}
	\pid(\thetav,\sigma^{2}\mid\dv) \propto L(\dv\mid\thetav,\sigma^{2})\pi_{0}(\thetav,\sigma^{2}).
\label{eqn:posterior}
\end{equation}
If we integrate out the dependence on $\sigma^{2}$, we are left with our target density
\begin{equation}
	\pi(\thetav\mid\dv) \propto L^{\sigma}(\dv\mid\thetav) \pi_{0}^{\lambda}(\theta_{1}) \pi_{0}^{I}(\theta_{2}),
\end{equation}
for $\thetav$ given $\dv$.
Exploiting the conjugacy of the inverse gamma and Gaussian distributions, we find that
\begin{equation}
	L^{\sigma}(\dv\mid\thetav) = t_{2\alpha_{\sigma}}\left(\dv \ \bigg| \ \calvG(\thetav), \frac{\beta_{\sigma}}{\alpha_{\sigma}}I_{\nd}\right)
\label{eqn:Lsigma}
\end{equation}	
is the density function of a multivariate $t$-distribution with $2\alpha_{\sigma}$ degrees of freedom, mean vector $\calvG(\thetav)$ and shape matrix $(\beta_{\sigma}/\alpha_{\sigma})I_{\nd}$. 

Determining the normalization constant
\begin{equation}
	Z:=\int_{\RR^{2}} L^{\sigma}(\dv\mid\thetav)\pi_{0}^{\lambda}(\theta_{1})\pi_{0}^{I}(\theta_{2})\dthetav,
\end{equation}
for $\pi$ gives the target density explicitly.
Determination of any QoI involving $\lambda$ and $I$ (expectations, variances, probabilities etc.), requires the computation of some integral involving the density $\pi$.
This task is non-trivial.
One approach, which is particularly amenable to situations with a small number of unknowns and where the target density is a smooth function of the unknowns is Smolyak sparse grid quadrature \cite{schillings2013sparse}.
An alternative approach, which we take here, is Markov chain Monte Carlo (MCMC) sampling  \cite{brooks2011handbook}.
MCMC methods allow us to sample from a distribution whose density is known only up to a constant of proportionality.
Our specific approach is described in Section \ref{subsubsec:mcmc}.

\subsubsection{Approximating the target density.}
\label{subsubsec:approximatingthetargetdensity}

Since we are unable to solve \eqref{eqn:pde}--\eqref{eqn:bc0} exactly, we are unable to evaluate $u$ in (\ref{eqn:ubar}) for a given value of $\thetav$.
When using an  MCMC method, the standard approach is to approximate the forward operator on a `case-by-case' basis.
That is, for each proposed $\thetav$, a new spatio-temporal approximation is computed.
We will refer to this as \textit{plain MCMC}, as in \cite{hoang2013complexity}.
For time-dependent and nonlinear problems, each forward solve usually requires the solution of a \textit{sequence} of linear systems.
This exacerbates the problem, and the cost of repeated forward solves can rapidly exceed a user's computational budget.

We propose an alternative approach where a surrogate is used to replace the repeated forward solves for different values of $\thetav$.
In particular, we use a pre-computed SGFEM approximation $\uhkt$ to the parametric forward problem \eqref{eqn:ypde}--\eqref{eqn:yq}.
As explained in Section \ref{subsec:forwardproblem}, once computed, this can be evaluated for any choice of the parameters $y_{1}$ and $y_{2}$ (corresponding to proposed values of $\lambda=\exp(\theta_{1})$ and $I=\exp(\theta_{2})$) without the need for any further linear system solves.
We will refer to this approach as \textit{SGFEM MCMC}. 

In detail, we use $\uhkt$ to approximate $L^{\sigma}(\dv|\thetav)$ in (\ref{eqn:Lsigma}) for proposed values of $\thetav$.
However, recall that $\uhkt$ is a function of $\yv\in\Gamma$.
From \eqref{eqn:ylambdaI}, each value of $y_{1}$ and $y_{2}$ can be mapped to values of $\lambda, I$ in the set $\Omega_{\lambda}\times\Omega_{I} := (\mu_{\lambda}-\nu_{\lambda}\sqrt{3}, \mu_{\lambda}+\nu_{\lambda}\sqrt{3}) \times (\mu_{I}-\nu_{I}\sqrt{3}, \mu_{I}+\nu_{I}\sqrt{3})$ via the mapping
\begin{equation}
	\zeta(\yv) := \left( \mu_{\lambda} + \nu_{\lambda} y_{1}, \mu_{I} + \nu_{I} y_{2} \right),
\label{eqn:zetamap}
\end{equation}
which has inverse
\begin{equation}
	\zeta^{-1}(\lambda,I) := \left( (\lambda-\mu_{\lambda})/\nu_{\lambda}, (I-\mu_{I})/\nu_{\lambda} \right).
\label{eqn:zetainvmap}
\end{equation}
Notice that $\yv\in\Gamma$ if, and only if, $\zeta(\yv) \in \Omega_{\lambda}\times\Omega_{I}$.

For values of $\thetav$ for which $(\lambda,I) \in \Omega_{\lambda}\times\Omega_{I}$, we replace $\calvG$ in (\ref{eqn:Lsigma}) with the \textit{approximate} observation operator given by
\begin{equation}
	\calvGhkt(\thetav)
	:= \big(\bar{u}_{hk\tau}(t_{1}, \zeta^{-1}(e^{\thetav})), \dots, \bar{u}_{hk\tau}(t_{\nd},\zeta^{-1}(e^{\thetav}))\big)^{\top},
\label{eqn:Ghk}
\end{equation}
where
\begin{equation}
\bar{u}_{hk\tau}(t,\yv):=\frac{1}{|D_{L}|}\int_{D_{L}} \uhkt(\rv, t, \yv) \, \mbox{d}S_{z}(\rv).
\label{eqn:uhktbar}
\end{equation}
This induces an approximate target density \cite{cotter2010approximation,dashti2016bayesian,stuart2010inverse} given by
\begin{equation}
\pihkt(\thetav\mid\dv) \propto L_{hk\tau}^{\sigma}(\dv\mid\thetav) \pi_{0}^{\lambda}(\theta_{1}) \pi_{0}^{I}(\theta_{2})
\label{eqn:pihktz}
\end{equation}
where $L_{hk\tau}^{\sigma}$ is given by \eqref{eqn:Lsigma} with $\calvG$ replaced by $\calvGhkt$.
Note that we use a number of time steps $\nt$ which is a multiple of $(\nd-1)$ so that the measurement times $t_{n}$ are a subset of the time steps $\tau_{n}$ and we need not interpolate in time to compute $\bar{u}_{hk\tau}(t_{n}, \zeta^{-1}(e^{\thetav})).$

We have no guarantee that the surrogate is accurate outside of $\Omega_{\lambda}\times\Omega_{I}$. 
Hence, for values of $\thetav$ for which $(\lambda, I) \notin\Omega_{\lambda}\times\Omega_{I}$ (equivalently $\yv\notin\Gamma$),
we would instead need to approximate \mbox{$L^{\sigma}(\dv \mid \thetav)$} using a standard deterministic spatio-temporal discretization (as in the plain MCMC approach).
In that case, we denote the induced approximate target density (defined analogously to $\pihkt$ in \eqref{eqn:pihktz}) by $\piht$.  

To ensure that we can evaluate $\uhkt$ for as many values of $\lambda$ and $I$ proposed by the selected MCMC method as possible,
we need to choose $\mu_{\lambda}$, $\nu_{\lambda}$, $\mu_{I}$ and $\nu_{I}$ in \eqref{eqn:ylambdaI} appropriately.
That is, the distributions used for $\lambda$ and $I$ in the construction of the surrogate should sufficiently cover the supports of the prior distributions selected for the Bayesian inference.
The uniform distribution selected for $\lambda$ is shown in Figure \ref{fig:figure2b}.
Should the MCMC propose an appreciable number of samples outside of the region in which the surrogate can be used, leading to a large number of significantly more expensive likelihood approximations (as in the plain MCMC approach), then it may be necessary to compute a new surrogate approximation which is accurate on a larger region of the state space.
This was not necessary in any of the examples that we present in Section \ref{sec:results}.

\subsubsection{Markov chain Monte Carlo.}
\label{subsubsec:mcmc}

As mentioned in Section \ref{subsubsec:approximatingthetargetdensity}, we use an MCMC algorithm to sample from the approximate target distribution, with density $\pihkt$.
For clarity of exposition, we use the simple Random Walk Metropolis Hastings (RWMH) \cite{brooks2011handbook} algorithm, although more advanced MCMC algorithms may be used, for example MALA \cite{roberts1998optimal} or HMC \cite{girolami2011riemann}.

The RWMH algorithm proposes values of $\thetav$ from a proposal distribution $q$ which is symmetric around the current state (usually a Gaussian), and then accepts/rejects them in the so-called Metropolis step in such a way that the resulting samples are the states of a Markov chain with stationary density equal to the target density.
This means that after a sufficient number of samples have been produced, say $n_{B}$, states can be considered draws from the target distribution.
Further post-processing of the chain, such as \textit{thinning} (where only a subset of the states are retained), can be implemented in order to reduce the correlation between samples.
Pseudocode for the algorithm used is given in Algorithm 1.

\begin{algorithm}
	\caption{SGFEM (RWMH) MCMC algorithm to sample from the approximate target density $\pihkt$.}
	
	compute SGFEM solution $\uhkt$ (offline)
	
	draw initial state $\thetav^{(0)}\sim\pi^{\lambda}_{0}\times\pi^{I}_{0}$
	
	\For{$ m = 0,1,2,\dots, M-1$}
	{
		draw proposal $\thetav^{*} \sim q=\calN(\thetav^{(m)}, \beta^{2}I_{2})$
		
		\eIf{$\exp(\thetav^{*})\in\Omega_{\lambda}\times\Omega_{I}$}
		{
			set $\hat{\pi}_{m+1} := \pihkt\left(\thetav^{*}\mid\dv\right)$
		}
		{
			set $\hat{\pi}_{m+1} := \piht\left(\thetav^{*}\mid\dv\right)$
		}
		
		set $p := 1 \wedge \frac{\hat{\pi}_{m+1}}{\hat{\pi}_{m}} \cdot \frac{q\left( \thetav^{(m)}\mid\thetav^{*}\right)}{q\left(\thetav^{*}\mid\thetav^{(m)}\right)}$
		
		draw $u\sim\calU(0,1)$
		
		\eIf{$u<p$}
		{
			set $\thetav^{(m+1)} := \thetav^{*}$ (accept)
		}
		{
			set $\thetav^{(m+1)} := \thetav^{(m)}$ (reject)
		}
	}
	
	discard $\thetav_{0},\dots,\thetav_{\nB-1}$ (burn in)
	
	thin chain
	
	compute Monte Carlo approximations $\EE[\varphi(\lambda,I)] \approx \frac{1}{(M+1-n_{B})}\sum\limits_{i=n_{B}}^{M}\varphi\left(e^{\theta_{1}^{(i)}},e^{\theta_{2}^{(i)}}\right)$
	
	\label{alg:rwmh}
\end{algorithm}

The SGFEM solve is performed \textit{once} outside the Monte Carlo loop.
There are $n_{h}n_{k}$ equations to solve per time step and the cost of each of these solves is $\calO(\nh\nk)$, provided an optimal solver (see \cite{benner2015low,powell2009block}) is available. This is an offline computation.
Recall that $n_{k}$ is defined in (\ref{nk_def}). Since $k$ denotes the chosen polynomial degree (typically $k=6$ is sufficiently high in this work), then $n_{k}$ is orders of magnitude smaller than $n_{h}$ and $n_{t}$.
The evaluation of the approximate observation operator $\calvGhkt$ is made once per iteration (online) and has only a small cost.
Specifically, since the observed quantity is averaged in space, evaluating $\bar{u}_{hk\tau}$ in (\ref{eqn:uhktbar}) for all $n_{z}$ measurement times only requires a matrix-vector product with a pre-computed matrix of size $n_{z} \times n_{k}$ with entries depending on the coefficients $u_{ij}^{n}$.
In each iteration, we simply need to compute a vector of length $n_{k}$ containing the Legendre polynomials evaluated at the point $\yv \in \Gamma$ corresponding to the proposed $\thetav$.
The total cost of the SGFEM MCMC approach consists of the offline cost (building the surrogate) and the online cost (evaluating the surrogate).
This is summarized in Table \ref{tab:table1}.

In the plain MCMC approach, we can use the same finite element discretization and time-stepping method as in the SGFEM approach.
At each iteration, however, we need to perform a forward solve to compute a new spatio-temporal approximation and evaluate it at the proposed value of $\thetav$.
The cost of each solve behaves as $\nt\times\calO(\nh)$, assuming an optimal solver for the FEM linear systems is available.
Note that we do not include the cost of evaluating the approximation in Table \ref{tab:table1} as the solve is the dominant cost.

It is clear that so long as the number of MCMC samples is larger than $\nk$ (the dimension of the polynomial space used to construct the surrogate), then a saving will be made with the SGFEM MCMC approach.
Due to the slow convergence of MCMC methods, $M$ is typically of order $10^{6}$ or higher whereas, since we have only two unknowns, $\nk$ only grows like $\calO(k^{2})$, where $k$ is the chosen polynomial degree.
Furthermore, as we save computational effort on every iteration, the larger $M$ is, the greater the savings made.

\begin{table}
	\centering
	\caption{\label{tab:table1} Computation costs of plain and SGFEM MCMC.}
	\begin{tabular}{@{}lcc}
	\hline
	&	 Offline 		&  Online  \\
	\hline
	plain  	&	  -/- 				&  $M\times \nt \times\calO(\nh)$ \\
	SGFEM 	& $\nt \times \calO(\nh \nk)$	&  $M\times\calO(\nd\nk)$ \\		
	\hline
	\end{tabular}
\end{table}

\subsubsection{Convergence and error.}
\label{subsubsec:convergence}

The computational cost of the SGFEM approach, but also the accuracy, depends on our choice of discretization parameters.
Since we use piecewise linear finite elements, the spatial approximation error for the forward problem can be expected to behave as $\calO(h)$, where $h$ is the mesh size.
The implicit Euler method yields an error which behaves as $\calO(\tau)$, where $\tau$ is the time step size. For parabolic PDEs, under certain simplifying assumptions, the error associated with the polynomial approximation on the parameter domain decays exponentially with respect to the polynomial degree $k$ (see \cite{nobile2009analysis,xiu2009efficient}).
Finally, the Monte Carlo error is $\calO(M^{-1/2})$, where $M$ is the number of samples.

In Section \ref{sec:results}, we demonstrate that the approximation to the target density $\pihkt$ obtained by the SGFEM MCMC approach converges as the polynomial degree $k$ is increased.
In future work, using the framework in \cite{cotter2010approximation} and \cite{hoang2013complexity} we hope to show convergence in the Hellinger distance \cite{stuart2010inverse} of $\pihkt$ to the true target density $\pi$ as the forward approximation is improved by decreasing $h$ and $\tau$ and increasing $k$.
In fact, under certain conditions \cite{cotter2010approximation} on the approximation of $\calvG$ by $\calvGhkt$, the rate of convergence in the forward approximation may be preserved in that of the target density.
This means that we are able to sample from a distribution with density close to $\pi$ as long as the forward approximation is sufficiently accurate.

\section{Results}
\label{sec:results}

A laser flash experiment was performed and temperature measurements were taken at $\nd=401$ equally spaced time points for the duration $T=40$ milliseconds.
The chosen values of all the model inputs are given in Table \ref{tab:table2}.
Prior distributions are stated for the unknowns $\lambda$ and $I$ rather than fixed values.

\begin{table}
	\centering
	\caption{\label{tab:table2} Parameter values for a single layer copper sample laser flash experiment.}
	\begin{tabular}{@{}lcc}
	\hline
	Quantity\,/\,units									& Symbol	& Value \\
	\hline
	Sample Radius\,/\,m 								& $R$ 		& $1.240\times10^{-2}$ \\		
	Sample Height\,/\,m 								& $H$ 		& $2.037\times10^{-3}$ \\
	Furnace/Ambient Temperature\,/\,K 					& $\Ta$ 	& $3.850\times10^{+2}$ \\
	Experiment Duration\,/\,s 							& $T$ 		& $4.000\times10^{-2}$ \\
	Laser Flash Duration\,/\,s	 						& $\tf$ 	& $4.000\times10^{-4}$ \\
	Laser Flash Depth\,/\,m 	 						& $\zf$ 	& $1.273\times10^{-4}$ \\
	Density\,/\,kg\,m$^{-3}$ 							& $\varrho$ & $8.930\times10^{+3}$ \\
	Specific Heat Capacity\,/\,J\,kg$^{-1}$\,K$^{-1}$ 	& $\cp$ 	& $3.970\times10^{+2}$ \\
	Heat Transfer Coeff.\,/\,W\,m$^{-2}$\,K$^{-1}$		& $\kappa$ 	& $1.100\times10^{+3}$ \\
	\hline
	Thermal Conductivity\,/\,W\,m$^{-1}$\,K$^{-1}$ 		& $\lambda$ & $\mbox{logN}(m_{\lambda},s_{\lambda})$ \\
	Laser Intensity\,/\,W\,m$^{-3}$ 					& $I$ 		& $\calU(0,\infty)$ \\
	\hline
	\end{tabular}
\end{table}

First, we use the SGFEM MCMC scheme to sample from the approximate target density $\pihkt$ in the case of a uniform laser profile.
That is, when $\chi(r)=1$ in the definition of the source term $Q$ in (\ref{eqn:q}) and (\ref{eqn:yq}).
We comment on the speedup achieved over the plain MCMC approach and demonstrate numerical convergence with respect to the polynomial degree $k$.
Second, we look at how the modelling of the laser profile affects the results of the Bayesian inference, comparing the results for the uniform profile to those obtained with two Gaussian-shaped profiles.

\subsection{Sampling the approximate target distribution using SGFEM MCMC}
\label{subsec:results1}

To construct the surrogate $\uhkt$, we use a spatial mesh with $\nh=12,206$ vertices and characteristic element size $h=4.24\times10^{-5}$\,m.
The number of implicit Euler time steps is set to $\nt=800$.
After investigating the accuracy of the solution to the forward problem, the polynomial degree for the parametric approximation is chosen to be $k=6$.
The resulting polynomial space has dimension $\nk=28$.
The time taken\footnote{CPU time. All experiments run on MacBook Pro laptop with 2.5\,GHz Intel Core i7 processor and 16\,GB RAM.} to complete the offline stage was 378 seconds.
 
\begin{figure*}
	\includegraphics[width=\textwidth]{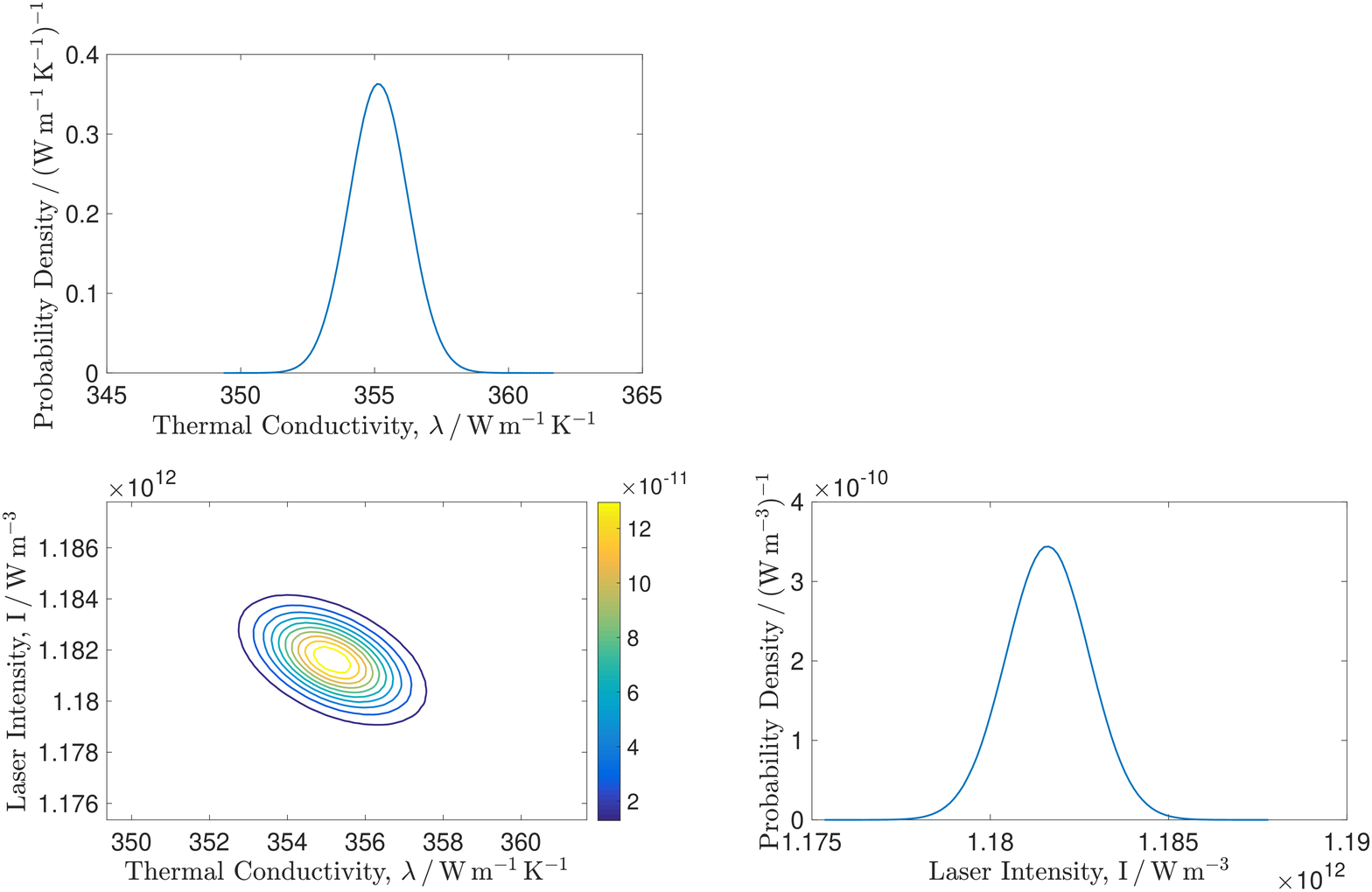}
	\caption{Approximate joint target density $\pihkt(\lambda,I \mid \dv)$ (bottom left) and approximate marginal densities $\pihkt^{\lambda}(\lambda \mid \dv)$ (top left) and $\pihkt^{I}(I \mid \dv)$ (bottom right) produced using 100 million MCMC samples with SGFEM surrogate computed using $n_{h}=12,206$, $n_{t}=800$ and $k=6$. Uniform laser profile $\chi(r)=1$. }
\label{fig:figure3}
\end{figure*}

The surrogate RWMH algorithm described in Algorithm 1 in Section \ref{subsubsec:mcmc} was run to generate 100 million samples from the approximate target distribution.
The time taken to produce these samples was 13,223 seconds.
We discarded the initial 10,000 samples as burn-in and the acceptance rate was tuned to near the optimal value of $23\%$ (see \cite{roberts2001optimal}) by choosing the proposal standard deviation $\beta=1.55$.
Histograms of the 100 million samples produced by Algorithm 1 were formed to obtain the approximate target density $\pihkt(\lambda, I \mid \dv)$ and the approximate marginal densities $\pihkt^{\lambda}(\lambda \mid \dv)$ and $ \pihkt^{I}(I \mid\dv)$ for $\lambda$ and $I$, respectively.
These are shown in Figure \ref{fig:figure3}.

Combining the offline and online steps, the total time required to approximate the target distribution was 13,601 seconds ($1.36\times10^{-4}$ seconds per sample).
For a finer spatio-temporal discretization, the cost of the offline stage would increase but, crucially, the cost of the online stage would remain unchanged (see Table \ref{tab:table1}).
Like all MCMC algorithms, this approach is highly parallelizable.
Once $\uhkt$ is computed, we could use it to produce multiple MCMC chains in parallel and combine the resulting samples appropriately.
However, we did not do this here.

No experiment was performed with the plain MCMC method.
Since a single forward solve costs around 35 seconds when we use the same spatio-temporal discretization, we estimate that the (CPU) time required to compute 100 million samples (again without exploiting parallelization) would be around 111 years!
This clearly shows that the surrogate approach represents a huge time saving (here, a reduction of 5 orders of magnitude) in comparison to the plain MCMC approach.
Even with a more modest 100,000 samples, the plain MCMC approach would take 40.5 days, compared to 391 seconds with the SGFEM approach.
With this smaller number of samples, the Monte Carlo error would far exceed the difference between the two approximate target densities $\piht$ and $\pihkt$.

The distributions shown in Figure \ref{fig:figure3} now tell us about the values of the unknowns $\lambda$ and $I$ (conditional on the data) and the uncertainty in these values.
We also gain information about the relationship between  $\lambda$ and $I$.
Note that this is a huge advantage over optimization approaches.
We obtain a mean value of $355.15$\,W\,m$^{-1}$\,K$^{-1}$  for $\lambda$ and a mean value of $1.1816\times10^{12}$\,W\,m$^{-3}$ for $I$.
The standard deviations are $1.10$ W~m$^{-1}$~K$^{-1}$ and $1\times10^{9}$ W~m$^{-3}$, respectively.  The mean value for $\lambda$ lies within the support of the prior (see Figure \ref{fig:figure2b}), but is different to the prior mean $\mu_{\lambda}=328.5$ W~m$^{-1}$~K$^{-1}$.
Moreover, we observe that the standard deviation is much reduced from the prior value $\sigma_{\lambda}=50$ W~m$^{-1}$~K$^{-1}$.

\begin{figure}
	\centering
	\includegraphics[width=\textwidth]{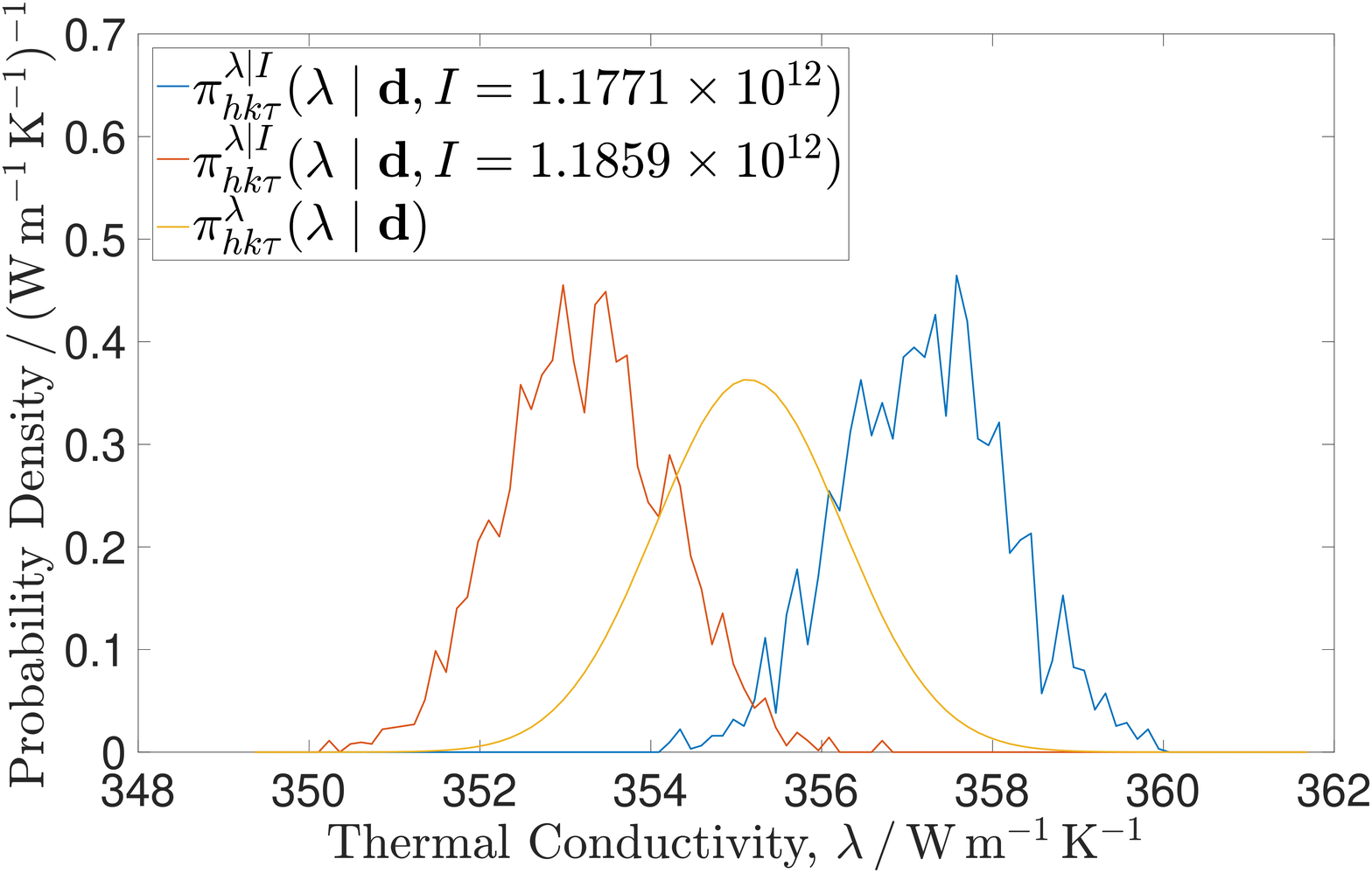}
	\caption{Approximate conditional densities \mbox{$\pihkt^{\lambda\mid I}(\lambda\mid\dv,I)$} for two different values of $I$ with the approximate marginal density $\pihkt^{\lambda}(\lambda \mid \dv)$ from Figure \ref{fig:figure3}. Uniform laser profile $\chi(r)=1$.}
	\label{fig:figure4}
\end{figure}

From Figure \ref{fig:figure3} we see that there is negative correlation between $\lambda$ and $I$ in the target distribution.
That is, a larger value of the thermal conductivity $\lambda$ requires a smaller value of the laser intensity $I$ (and vice versa) for the model output to be as consistent with the data as possible.
The computed value of the correlation coefficient is -0.48.
This highlights the importance of performing inference on both unknowns simultaneously, as fixing the value of one amounts to considering the conditional distribution for the other.
The approximate densities of two such conditional distributions are displayed in Figure \ref{fig:figure4}, along with the approximate marginal target density for $\lambda$.
Observe that the three densities indicate very different likely values of $\lambda$. Note also that fixing a lower (resp. higher) value of the laser intensity $I$ results in a larger (resp. smaller) likely value of $\lambda$ being indicated. 
The histograms used to approximate the conditional densities are far less converged than those for the marginals as there are fewer samples which may be used to construct them.
Conditional distributions for larger or smaller values of $I$ than those used here did not have sufficient samples to create histograms, but would result in more severe differences in the distributions on $\lambda$. Correlation between $\lambda$ and $I$, and hence the need for joint inference, has not previously been considered for the data set considered here.

\begin{figure*}
	\includegraphics[width=\textwidth]{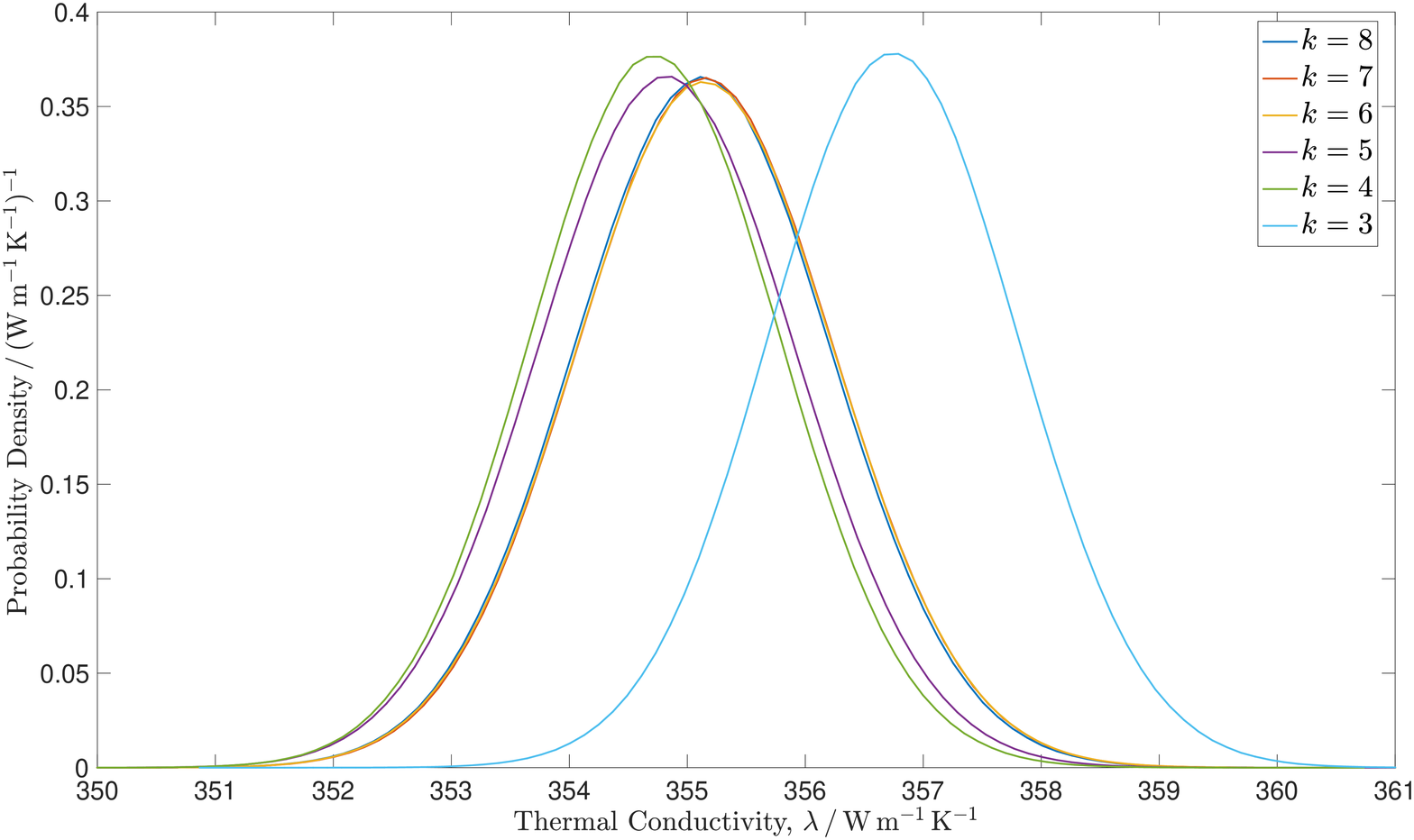}
	\caption{Approximate marginal target densities $\pihkt^{\lambda}(\lambda \mid \dv)$ for $\lambda$ produced using 100 million MCMC samples with SGFEM surrogate computed using $n_{h}=12,206$, $n_{t}=800$ and varying values of $k$. Uniform laser profile $\chi(r)=1$.}
	\label{fig:figure5}
\end{figure*}

In Figure \ref{fig:figure5}, we show the results of a numerical convergence study.
We plot the approximate marginal target density for $\lambda$ obtained using a varying value of $k$ in the construction of the surrogate.
The spatio-temporal discretization is fixed as before.
For $k \ge 6$, the densities produced are almost indistinguishable visually.
The approximation converges as $k$ increases.
Similar behaviour is witnessed when the discretization parameters $h$ and $\tau$ are decreased.

\begin{figure*}
	\includegraphics[width=\textwidth]{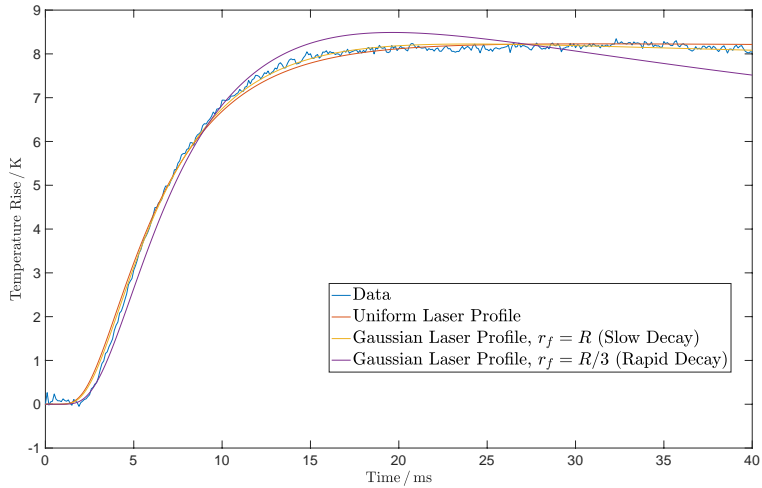}
	\caption{Thermogram of experimental data $\dv$ and model thermograms for three different laser profiles computed by evaluating the observation operator $\calvGhkt$ at the approximate target distribution mean obtained using 100 million MCMC samples and SGFEM surrogate computed using $n_{h}=12,206$, $n_{t}=800$ and $k=6$.}
	\label{fig:figure6}
\end{figure*}

In Figure \ref{fig:figure6}, we compare the experimental data $\dv$ with a model thermogram computed by evaluating the approximate observation operator $\calvGhkt$ in \eqref{eqn:Ghk} at the approximate target distribution mean (computed with $k=6$).
We observe that, modelling the laser profile as uniform (i.e., $\chi(r) = 1$), a very good fit to the data has been achieved.
This could potentially be improved further by treating other model inputs as unknown (such as the heat transfer coefficient $\kappa$ in the boundary condition) and inferring their values from the data too.

\subsection{Varying the laser profile}
\label{subsec:results2}

In this section, we consider the effect on our inference results of using a Gaussian-shaped profile to describe the laser flash, instead of a uniform profile.
That is, we now set
\begin{equation}
	\chi(r) := \exp(-r^{2}/2\rf^{2})
\label{eqn:chiGaussian}
\end{equation}
in \eqref{eqn:yq}, rather than $\chi(r) =1$ as in Section \ref{subsec:results1}.

By varying the parameter $\rf$ in \eqref{eqn:chiGaussian}, we are able to control how quickly the laser intensity in the model decays as the distance from the centre of the sample increases.
We repeat the experiment in Section \ref{subsec:results1}, now using the values $\rf=R$ and $\rf=R/3$ to represent lasers whose profiles decay slowly and rapidly, respectively.
The resulting estimates of the marginal target distribution means and standard deviations for $\lambda$ are displayed in Table \ref{tab:table3}. 
Note that such investigations would be impossible without the rapid sampling provided by the SGFEM surrogate approach.
We see that the profile of the laser has a significant effect on the reported values.
The mean value of $\lambda$ decreases as the laser becomes more focused at the centre of the sample.
This indicates that checking the uniformity of the laser is important for accurate measurement and that uncertainties about the laser profile shape should be built into the model.
Note that the computed values for $I$ are not directly comparable. 

\begin{table}
	\centering
	\caption{\label{tab:table3} Approximate marginal target distribution means and standard deviations for $\lambda$ for different laser profiles.}
	\begin{tabular}{@{}lccc}
	\hline
	 & Uniform 					& $\rf = R$ & $\rf = R/3$  \\
	\hline
	\ \ \, mean $\lambda$\,/\,W\,m$^{-1}$\,K$^{-1}$ 	   & $355.15$ 			& $336.86$					& $281.94$ \\
	st. dev. $\lambda$\,/\,W\,m$^{-1}$\,K$^{-1}$  & $1.10$  			& $0.74$					& $1.61$ \\
	\hline
	\end{tabular}
\end{table}

Figure \ref{fig:figure6} shows model thermograms computed by evaluating $\calvGhkt$ at the target distribution mean for all three laser profiles.
It is clear that a good level of fit to the experimental data has been obtained when using the uniform laser profile or the flatter Gaussian profile ($\rf=R$).
When using the rapidly decaying Gaussian profile however, the target distribution mean gives a bad fit to the data.
This indicates that modelling the laser this way is incorrect and that the experiment has been performed in such a way that the laser profile is actually close to the ideal uniform shape.

\begin{figure}
	\centering
	\includegraphics[width=0.75\textwidth]{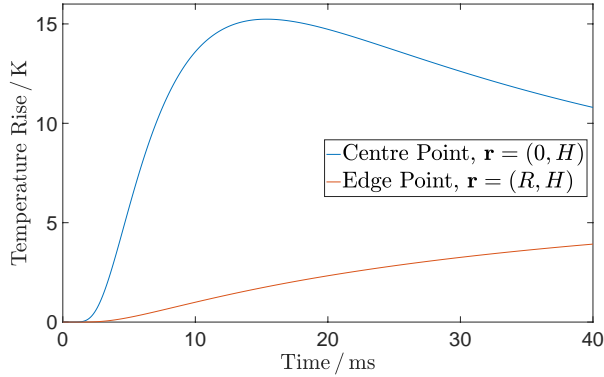}
	\caption{Approximate temperature $\uhkt$ evaluated at the target density mean for $(\lambda,I)$ at $\rv=(0,H)$ (centre of the top face) and at $\rv=(R,H)$ (outer edge of top face) for the Gaussian laser profile with fast decay, i.e., $\chi(r):=\exp(-9r^{2}/2R^{2}$).}
\label{fig:figure7}
\end{figure}

The bad fit observed for the sharply decaying Gaussian profile is indicative of the fact that in this model the heat energy travels rapidly through the centre of the sample and then slowly spreads out in the radial direction.
The simulated temperature does not reach a steady state over the time scale of the experiment.
Figure \ref{fig:figure7} shows the approximate temperature at the spatial locations $(r,z)=(R,H)$ and $(r,z)=(0,H)$ obtained by evaluating $\uhkt$ at the target density mean for the duration of the experiment.
At the point $(r,z)=(0,H)$, the temperature increases rapidly initially and is then is still cooling at time $T$. However, at $(r,z)=(R,H)$, the temperature increases more slowly and is still increasing at time $T$.

Note that, for each laser profile considered, there is a discrepancy at the early times where the measured temperature spikes.
This is due to the fact that the laser flash itself is detected by the sensor.
This feature has not been incorporated into the model.

\section{Conclusions and future work}
\label{sec:conclusions}

In this paper we formulated the determination of the thermal diffusivity (and hence the conductivity) of a material given laser flash experiment data as a Bayesian inverse problem.
We demonstrated how an SGFEM surrogate may be used to accelerate MCMC sampling from an approximate target distribution, in a real life situation where the plain MCMC approach is infeasible.
We explained how posing the problem as a Bayesian inverse problem allows for clear quantification of the uncertainty in the estimate of the thermal conductivity.
We used this approach to find the joint distribution on the thermal conductivity and the laser intensity.
We observed a correlation in the posterior between these unknowns, which indicates that incorrect estimates of the laser intensity without joint inference can lead to biased estimates of the thermal conductivity.
A numerical convergence study was presented, showing the effect of varying the discretization parameters in the construction of the surrogate.

Due to the speedup over plain MCMC provided by the SGFEM MCMC approach, we have been able to investigate the effect of the laser profile on the inference results.
In particular, we have shown that the spatial uniformity of the laser profile strongly affects the value of the estimated thermal conductivity.
This demonstrates the importance of both accurately modelling the laser flash and correctly setting up the physical experiment.

Combining a parametric surrogate with MCMC sampling offers a computationally efficient way to solve inverse problems involving PDEs.
However, we caution that the appropriate choice of surrogate is problem-dependent.
For problems where the forward solution is a non-smooth function of the parameters chosen to represent the uncertain inputs, the standard SGFEM approach adopted here is not recommended. Moreover, for problems where the solution becomes more nonlinear over time, using a fixed polynomial basis for all time steps, as has been done here, is not a sensible approach. When a surrogate is chosen, its accuracy with respect to the solution of the forward problem should first be investigated before using it to accelerate the numerical solution of an inverse problem.

In future, we shall consider how the multi-thermogram approach of \cite{allard2015multi} can similarly be incorporated into the SGFEM MCMC framework (for single and multi-layered materials) through simultaneous evaluation of surrogates for multiple forward problems.
We will also incorporate more uncertain inputs into the model, for example the heat transfer coefficient $\kappa$ in the definition of the boundary condition.

\section*{Acknowledgements}
The authors would like to acknowledge the support of the National Measurement System of the UK
Department of Business, Energy \& Industrial Strategy, which funded this work as part of NPL's Data Science program.

\bibliographystyle{siam}
\bibliography{myBibliography}

\end{document}